\magnification\magstep1

\centerline{\bf The Number of Periodic Orbits of a }

\centerline{\bf Rational Difference Equation}
\bigskip
\centerline{ Eric Bedford\footnote*{Supported in part by the NSF.}  and Kyounghee Kim}
\bigskip
\centerline{\it Dedicated to Christer Kiselman}

\medskip
\bigskip
\noindent{\bf \S0.  Introduction } 
Here we consider difference equations of the form
$$z_{n+1}=g(z_n,z_{n-1}) = { p(z_n,z_{n-1})\over q(z_n,z_{n-1})}\eqno(0.1)$$
where $g$ is a rational function of two variables, which we write as a quotient of polynomials without common factor.  Given starting values $z_0=x$ and $z_1=y$, equation (0.1) gives rise to an infinite sequence $(z_n)_{n\ge0}$, as long as the denominator in (0.1) does not vanish.   One of the basic questions here is to find the periodic sequences generated by (0.1), i.e., sequences for which $z_{j+p}=z_j$ for all $j$.  We refer to the books [KL], [KM] and [GL] for expanded discussions of this question.  Here we describe a method for giving an upper bound on the complex periodic sequences of period $p$.  We note that equations (0.1) have been widely studied in the real domain, frequently for $z_n>0$;    when we obtain an upper bound on the number of complex periodic points, this upper bound  also applies to the periodic sequences which are real (or positive).

Let us write (0.1) as a rational map
$$F(x,y)=(y,g(x,y))\eqno(0.2)$$
of the plane.   It is most natural for us to consider $F$ as a map of the complex plane ${\bf C}^2$.  In fact, we will work with the extension of $F_X:X\to X$ to a compactification $X$ of ${\bf C}^2$.  We let $Ind(F)$ denote the set of points of indeterminacy of $F$.  The set $Ind(F)$ is finite and has the properties: (1) $F$ is holomorphic on $X-Ind(F)$, and (2) $F$ cannot be extended  continuously to any point of $Ind(F)$.  (In fact, $F$ blows $q$ up to a subvariety of $X$.)   We will consider the set
$${ Fix}'_n=\{z\in{\bf C}^2: F^jz\in{\bf C}^2- Ind(F)\ \forall\ 1\le j\le n, {\rm\ and\ }f^nz=z\}.$$
Thus a point $(x,y)$ belongs to $Fix'_n$ if and only if the corresponding sequence $(z_n)$ is periodic for some period $p$ which divides $n$.  The mappings $F$ that we consider here will have only isolated periodic points, and each such point $q$ has a well-defined (algebraic) multiplicity $\mu_q$, which at a regular point is the multiplicity of $q$ as a solution of the pair of equations $(x,y)-F^n(x,y)=(0,0)$.  We set 
$$\# { Fix}'_n  =  \sum_{q\in Fix'_n} \mu_q$$ 
denote the number of fixed points counted with multiplicity.  The multiplicity is $\ge1$, so an upper bound on $\#{ Fix}'_n$ gives an upper bound on the number of points in ${ Fix}'_n$.

We will  consider  linear fractional recurrences of the form
$$z_{n+1}={a_0 + a_1 z_{n-1} + a_2 z_n\over b_0 + b_1 z_{n-1} + b_2 z_n}.  \eqno(0.3)$$
For instance, the host-parasite model (see [KM, p.\ 181]), which is customarily presented as 
$$\eqalign{x_{n+1} & = {\alpha x_n\over 1+\beta y_n} \cr 
y_{n+1} &={\gamma x_n y_n \over 1+\beta y_n}} \eqno(0.4)$$
may be re-written in the form (0.3) (see \S2).
Our first result is:
\proclaim Theorem 1.  The  periodic solutions of (0.3) satisfy
$$\#{Fix}_n'\le F_{n+1}+F_{n-1},\eqno(0.5)$$
where $F_n$ denotes the $n$th Fibonacci number.  For generic  $a_j$ and $b_j$, equality holds in (0.5).

We may also consider the special case
$$z_{n+1}={a+z_n\over b+z_{n-1}}. \eqno(0.6)$$
Let us define a Fibonacci-like sequence defined by $\varphi_{n+3}=\varphi_n+\varphi_{n+1}$ and with starting values $\varphi_1=0$, $\varphi_2=2$, $\varphi_3=3$.
\proclaim Theorem 2.  The periodic points of (0.6) satisfy
$$\#{ Fix}'_n\le \varphi_{n}+2. \eqno(0.7)$$
Further, equality holds in (0.7) for generic $a$ and $b$.

The method of proof is an adaptation of the Lefschetz Fixed Point Theorem.  In \S1 we discuss this method, which was used for a different family of rational maps in [BD].  In \S2 we give the proofs of Theorems 1 and 2.  And in \S3 we show how it may be applied to two other examples.

\bigskip\noindent{\bf \S1.  Lefschetz Fixed Point Theorem }  Let us consider a polynomial map $F:{\bf C}^2\to{\bf C}^2$ which is given as
$$F(x,y)=(f_d(x,y),g_d(x,y))+\cdots,\eqno(1.1)$$
where $f_d$ and $g_d$ are homogeneous polynomials of degree $d$, and the dots indicate terms of lower degree.  We suppose that $\{f_d=g_d=0\}=\{(0,0)\}$.  Such a polynomial map has a continuation to a holomorphic map $F:{\bf P}^2\to {\bf P}^2$, where ${\bf P}^2$ denotes the projective plane, with points represented by homogeneous coordinates $[x_0:x_1:x_2]$.  Thus ${\bf P}^2={\bf C}^2\cup\Sigma_0$, where we identify points of ${\bf C}^2$ with points of ${\bf P}^2$ by the map $(x,y)\mapsto [1:x:y]$, and $\Sigma_0=\{x_0=0\}$ is the hyperplane at infinity.

The cohomology of ${\bf P}^2$ is $H^*({\bf P}^2;{\bf Z})\cong H^{0}\oplus H^2\oplus H^4\cong {\bf Z}\oplus {\bf Z}\oplus {\bf Z}$.  For any rational map, the action of $F^*$ on $H^0$ is the identity and $F^*$ acts by multiplication by the topological degree $d_{\rm top}$ on $H^4$.  Since $F$ is holomorphic and is homogeneous of degree $d$,  we have $d_{\rm top}=d^2$, and $F^*$ acts as multiplication by $d$ on $H^2$.  Thus the total map on cohomology is given by the diagonal matrix
$$F^*|H^*({\bf P}^2)=\pmatrix{1 &&\cr
&d &\cr
&& d^2}.$$
If $d>1$, then the fixed points of $F$ are isolated, and so the Lefschetz Fixed Point Theorem  tells us that the number of fixed points is given by the trace of the mapping on cohomology:
$$\# Fix(F)=Tr(F^*|H^*({\bf P}^2)) = 1+d+d^2.$$
Here the expression $\# Fix(F)$ denotes the sum of the multiplicites of the fixed points in ${\bf P}^2$.
If we wish to find $Fix'_1(F)$, then we consider which of the fixed points are at infinity.  We see that $\Sigma_0$ is invariant, and the restriction $F|\Sigma_0$ is given by the rational map $[0:x_1:x_2]\mapsto [0:f_d(x_1,x_2):g_d(x_1,x_2)]$.  If we count the fixed points of $F|\Sigma_0$ using the Lefschetz Fixed Point Theorem argument above, we find
$\#Fix(F|\Sigma_0)=1+d$.  Thus we can count the fixed points in ${\bf C}^2$ as
$\#Fix'_1(F)=d^2$.  

Now we observe that $F^n$ is a mapping satisfying (1.1), and the degree of $F^n$ is $d^n$.  Applying our formula, we have $\# Fix'_n=d^{2n}$.
\bigskip
Let us next consider a mapping that arises in the so-called $SI$ epidemic model (see [KM, page 186]):
$$g: (S,I)\mapsto (S+\alpha SI, I-\alpha SI).\eqno(1.2)$$
We note that a related model, of the form
$$(S,I)\mapsto (S+\alpha SI+\lambda I, (1-\lambda)I - \alpha SI)\eqno(1.3)$$
can be treated in the same way.
We may extend $g$ to projective space by setting $S=x_1/x_0$ and $I=x_2/x_0$.  This gives
$$g[x_0:x_1:x_2] = [x_0^2: x_1(x_0-\alpha x_2): x_2(x_0-\alpha x_1)].$$
We see that $g$ defines a holomorphic map of ${\bf P}^2$ except at the indeterminacy locus 
$$Ind(g)=\{e_1=[0:1:0],e_2=[0:0:1]\}.$$
We see that both points $e_1$ and $e_2$ blow up to the hyperplane at infinity $\Sigma_0$.  The hyperplane at infinity is exceptional in the sense that
$$\Sigma_0-\{e_1,e_2\}\mapsto [0:1:-1].$$
Further, we note that $[0:1:-1]\in\Sigma_0$ and is a fixed point of $g$.   

We observe that the algebraic degree of $g$ is 2, and the topological degree is also 2.  Thus we have
$$g^*|H^*({\bf P}^2) = \pmatrix{1 &&\cr
&g^*|H^{2}&\cr
&& {\rm deg}_{\rm top}\cr}= \pmatrix{1 &&\cr
&2&\cr
&& 2\cr}.$$
Now the trace of $g^*$ is equal to the intersection number between the diagonal $\Delta\subset{\bf P}^2\times{\bf P}^2$ and $g^{-1}\Delta$ (see [Fu]).  In the case of our mapping, the intersection points are isolated, and we have
$$Tr (g^*) = 1+2+2 = \sum_{q\in \Delta\cap g^{-1}\Delta}\mu_q.$$

The points of $\Delta\cap g^{-1}\Delta-Ind(g)$ are exactly the fixed points of $g$.  For a point $q\in( \Delta\cap g^{-1}\Delta)\cap Ind(g)$ there are two possibilities:  (1) $q$ blows up to a curve $\gamma\ni q$; or (2)  there is an exceptional curve $E$ which contains $q$ and $g(E-Ind(g))=q$; (or both).  We have seen above that $e_1$ and $e_2$ are in case (1):  they are indeterminate points which blow up to $\Sigma_0$, and $e_1,e_2\in\Sigma_0$.  We note that while $[0:1:-1]\notin Ind(g)$ it is in case (2), since it is the image of the exceptional curve $\Sigma_0$.  Further, since $\Sigma_0$ is mapped to itself, we see that there can be no orbit which intersects both ${\bf C}^2$ and $\Sigma_0$.   Thus there are no periodic points whose orbits pass through infinity.  This gives:
$$\# Fix'_1(g) = Tr(g^*|H^*)-3 = 2.$$
We observe that applying this argument to (1.2) and (1.3) is the same, since we are only concerned with the behavior of points at infinity, and this is not influenced by the lower degree terms of $g$.

In order to discuss $Fix'_n$, we consider the iterate $g^n$.  $\Sigma_0$ is the only exceptional curve for $g$, and $\Sigma_0$ is taken to a fixed point.  Thus there is no exceptional curve that is mapped to $Ind(g)$ by an iterate of $g$.  This is the criterion given in  [FS] for a map of ${\bf P}^2$ to satisfy  $deg(g^n) =(deg(g))^n=2^n$.   Thus we have
$$g^*|H^* = \pmatrix{1 &&\cr
& 2^n&\cr
&& 2^n\cr}.$$
and so the intersection number of $\Delta\cap g^{-n}\Delta$ is $1+ 2\cdot 2^n$.  Now the 3 points $e_1$, $e_2$ and $[0:1:-1]$ exhaust all the indeterminacy points and exceptional images, so there are no new ``fake periodic points'' for higher period.  This gives
$$\#Fix'_n(g) = 2\cdot 2^n-2.$$

\bigskip\noindent{\bf \S2.  Linear Fractional Recurrences }  We wish now to show how to apply the ideas of  \S1 to the maps (0.3) and (0.6).  Let us start by rewriting the recurrence (0.3) as a planar map of the form (0.2).  Rewriting this mapping in terms of homogeneous coordinates $[x_0:x_1:x_2] = [1:x:y]$, we have
$$F[x_0:x_1:x_2] = [x_0 \beta\cdot x: x_2\beta\cdot x:x_0\alpha\cdot x]$$
where $\alpha=(a_0,a_1,a_2)$ and $\beta=(b_0,b_1,b_2)$.  $F$ is invertible, so the topological degree is 1, and we have
$$F^*|H^*({\bf P}^2)=\pmatrix{1&&\cr
&2&\cr
&&1\cr}.$$
The points of indeterminacy are
$Ind(F) = \{e_1, p_0, p_\gamma\}$, where $p_0=[0:-\beta_2:\beta_1]$ and $p_\gamma=[\beta_1\alpha_2-\beta_2\alpha_1: -\beta_0\alpha_2+\alpha_0\beta_2:\alpha_1\beta_0-\alpha_0\beta_1]$.  
See [BK, \S1] for details.  For generic values of $\alpha$ and $\beta$, $e_1 \in\Delta\cap f^{-1}\Delta$ has multiplicity 2.  First, $e_1$ is a fixed point because $e_1\in\Sigma_0$, which is an exceptional curve, and $\Sigma_0\mapsto e_1$, and second it satisfies the property (2) above, because $e_1$ blows up to the curve $\Sigma_B$, and $e_1\in \Sigma_B$.  Taking the intersection number and subtracting the multiplicity of $e_1$, we have
$$\#Fix'_1(F)=Tr(F^*|H^*)-2= 2.$$
We note that for non-generic values of $\alpha$ and $\beta$ this becomes $\#Fix'_1(F)\le 2$.

We now want to treat higher iterates of $F$, but $deg(F^n)\ne(deg(F))^n$.   This is seen because, $\Sigma_0\mapsto e_1\in Ind(F)$, i.e.\ the forward orbit of an exceptional curve enters the indeterminacy locus.   We deal with this situation by using the approach of Diller and Favre [DF], in which we replace ${\bf P}^2$ by a space on which the passage to cohomology is consistent with iteration.  In this case, we work with the space $Y$ which is obtained by blowing up the point $e_1$; further details are given in [BK, \S2,3].   $H^{2}(Y;{\bf Z})$ is generated by the class of a hyperplane $H$ and the blowup fiber $E_1$.   With respect to this basis, the induced rational map $F_Y:Y\to Y$ acts on $H^2(Y;{\bf Z})$ according to
$$F_Y^*|H^{2}(Y)=\pmatrix{2 & 1\cr
-1 &-1\cr}.$$
The map $F_Y$ has two exceptional curves, which are denoted $\Sigma_\beta$ and $\Sigma_\gamma$, and $Ind(F_Y)$ consists of two points.  For generic $\alpha$ and $\beta$, the orbits of $\Sigma_\beta$ and $\Sigma_\gamma$ do not meet the indeterminacy locus, so the map $F_Y$ is regular on $H^{2}$, which means that $(F_Y^n)^*=(F_Y^*)^n$.  Using this, as well as the defining properties of the Fibonacci numbers $F_n$, we have
$$(F^n)^*=\pmatrix{2 & 1\cr
-1 &-1\cr}^n=\pmatrix{F_{n+2} &  F_n\cr
-F_n&-F_{n-2}\cr}.$$
Thus  we obtain
$$\#Fix'_n=F_{n+2}-F_{n-2}=F_{n+1}+F_{n-1}.$$
This proves Theorem 1.

\bigskip  {\it Another formulation.}  Let us show how to rewrite the mapping (0.4) as a (linear fractional) difference equation (0.3) with $a_0=a_1=0$, $a_2=\alpha$, $b_0=1$, $b_1=\gamma$, and $b_2=-\beta/\alpha$.  Applying  (0.4) we have
$$\gamma x_n-\beta y_{n+1}={\gamma x_n + \gamma\beta x_ny_n\over 1+\beta y_n} -{\gamma\beta x_n y_n\over 1+\beta y_n} = {\gamma x_n\over 1+\beta y_n} = {\gamma\over \alpha}x_{n+1}.$$
Thus 
$$y_{n+1} = {\gamma\over\beta}(x_n-{1\over \alpha}x_{n+1}),\eqno(2.1)$$
which means that the sequence $\{(x_n,y_n)\}$ in (0.4) is determined by the sequence $\{x_n\}$.  From the first equation in (0.4) and (2.1) we have
$$x_{n+1}= {\alpha x_n\over 1+\gamma x_{n-1}-{\beta\over \alpha} x_n},$$
so the sequence $\{x_n\}$ satisfies (0.3) with the given values of $a_j$ and $b_j$.

\bigskip {\it Proof of Theorem 2.}
Now let us consider the map (0.6), which is the mapping (0.3), but restricted to the (non-generic) case $\alpha=(a,0,1)$ and $\beta=(b,1,0)$.  There are three points of indeterminacy $Ind(F) = \{e_1,e_2,p\}$, and there are three exceptional curves $\Sigma_0$, $\Sigma_\beta$ and $\Sigma_\alpha$.  We have $\Sigma_\beta\mapsto e_2\in Ind(F)$ and $\Sigma_0\mapsto e_1\in Ind(F)$.   We let $X$ denote the space obtained from ${\bf P}^2$ by blowing up $e_1$ and $e_2$.  We let $F_X:X\to X$ denote the induced rational map.  $F_X$ has one indeterminate point $Ind(F_X)=\{p\}$.  The line $\Sigma_\alpha\subset X$ is exceptional, and $\Sigma_\alpha\mapsto q\mapsto F_Xq\mapsto F^2_Xq\mapsto\cdots$.   For generic values of $a$ and $b$, the orbit of the exceptional line does not land on $Ind(F_X)=\{p\}$.  In this case we have $(F_X^*)^n=(F_X^n)^*$.

If $H$ denotes the class of a general hyperplane in $X$, and if $E_1$ and $E_2$ are the exceptional blowup fibers, then $\{H,E_1,E_2\}$ gives an ordered basis for $H^2(X;{\bf Z})$.  In [BK, \S4] we saw that
$$F_X^*=\pmatrix{2 & 1& 1\cr
-1& -1 & 0\cr
-1& -1 & -1\cr}$$
with respect to this basis.

The characteristic polynomial of $A:=F_X^*$ is $x^3-x-1$ so we have $A^3-A-I=0$.  It follows that each entry $a_{i,j}^n$ of $A^n$ satisfies
$$a_{i,j}^{n+3}=a_{i,j}^{n+1}+a_{i,j}^n.$$
Now let $\varphi_n$ be the sequence satisfying the linear recurrence $\varphi_{n+3}=\varphi_{n+1}+\varphi_n$, and which has starting values  $\varphi_0=3\ (=Tr(I))$, $\varphi_1=0\ (=Tr(A))$, and $\varphi_2=2\ (=Tr(A^2))$.  It follows that $\varphi_n=Tr(A^n)$ for all $n$, which proves Theorem 2.

\bigskip\noindent{\bf \S3.  Two Examples }
In this section we will consider two maps that arise as biological models.  The first (see [KM, p.\ 174]) is
$$g:\ \ \ (x_{n+1},y_{n+1}) = \left(\alpha x_n + {x_n\over a_0 + a_1 x_n + a_2 y_n}, \beta y_n + {y_n\over b_0 + b_1 x_n + b_2 y_n}\right) $$
For this map, we determine the topological degree as well as the degree as a homogeneous polynomial.  Then we show that $deg(g^n)=(deg(g))^n$ and count the spurious fixed points. 

Finding the number of preimages of a generic point, we see that $d_{\rm top}=4$.  Checking for indeterminate points in ${\bf C}^2$, we see that we have fractions of the form ${0\over 0}$ at $p_a:=(0,-a_0/a_2)$ and $p_b:=(-b_0/b_1,0)$.   The point $p_a$ blows up to a horizontal line, which contains $e_1\in{\bf P}^2$; and $p_b$ blows up to a vertical line, which contains $e_2\in {\bf P}^2$.  

If we write $g$ in homogeneous coordinates, we find the mapping
$$g[x_0:x_1:x_2] = [x_0 \ell_a(x)\ell_b(x): x_1 \ell_b(x) m_1(x): x_2 \ell_a(x) m_2(x)]$$
where $\ell_a(x) = a_0x_0 + a_1 x_1 + a_2 x_2$,  $\ell_b$ is similar, and $m_1$ and $m_2$ are linear functions.  Thus $g$ has algebraic degree 3.   Let us set $L_a=\{\ell_a=0\}$ and $L_b=\{\ell_b=0\}$.  The points of indeterminacy are then $Ind(g)=\{p_a,p_b,q_a,q_b,r\}$, where $q_a=L_a\cap\Sigma_0=[0:-a_2:a_1]$, $q_b=L_b\cap\Sigma_0=[0:-b_2:b_1]$, and $r=L_a\cap L_b$.  We observe, also, that $p_a\in L_a$ and $p_b\in L_b$.   The  Jacobian determinant of $g$ is given by $\ell_a\ell_b R(x)$, where $R(x)$ is quartic.  In fact, $\{R=0\}$ is a curve of branch points.  Thus the set of exceptional curves consists of $L_a\mapsto e_1$ and $L_b\mapsto e_2$.  The points $e_1$ and $e_2$ are regular (not indeterminate) and are fixed.  Since the orbits of exceptional curves never enter the indeterminacy locus, it follows that $deg(g^n)=(deg(g))^n=3^n$.  Thus we have $g^{n*}|H^*({\bf P}^2)=\pmatrix{1&&\cr & 3^n&\cr && 4^n\cr}$.  

Now we count the points that are not in $Fix_n'$, but which will contribute to the intersection number of $\Delta\cap f^{-1}\Delta$.  Namely, we have the (regular) fixed points $e_1$ and $e_2$ at infinity.  Further, we have $p_a$, which blows up to $L_a\ni p_a$, and $p_b$, which blows up to $L_b\ni p_b$, for a total count of 4.   For generic parameters, we have $q_a\ne q_b$, etc, so taking the trace and subtracting the excess multiplicity, we have 
$$\#Fix_n'=3^n+4^n-3.$$

We apply a similar analysis to
$$f(x,y) = \left( {axy + x+ y\over b+x+y},{cxy + x+y\over d+x+y}\right)$$
(see [KM, p.\ 172]).  This map differs from the previous one because it needs to be regularized.

First, we note that the topological degree of $f$ is 2.
If we rewrite this mapping in terms of homogeneous coordinates $[x_0:x_1:x_2]=[1:x:y]$ on projective space ${\bf P}^2$ then we have
$$\eqalign{ f[x_0: & x_1:x_2] = [x_0(b x_0+x_1+x_2)(d x_0+x_1+x_2): \cr
&(dx_0+x_1+x_2)(ax_1x_2+x_0x_1 +x_0x_2): (bx_0+x_1+x_2)(cx_1x_2+x_0x_1 +x_0x_2)].\cr}$$
To find the exceptional curves, we look at the Jacobian, which is 
$$x_0(x_1-x_2)(b x_0+x_1+x_2)(d x_0+x_1+x_2)Q(x)$$
for some quadratic polynomial $Q(x)$.  The line $\{x_1=x_2\}$ consists of branch points;  all of the other critical curves are exceptional:
 $$\eqalign{ &\{x_0=0\}\mapsto [0:a:c]\in \Sigma_0, \ \ \{b x_0+x_1+x_2=0\}\mapsto e_1,\cr
& \{d x_0+x_1+x_2=0\}\mapsto e_2,\ \  \{Q(x)=0\}\mapsto [ac(b-d): a(a-c)d: b(a-c)c].\cr}$$
We observe that $[0:a:c]\in\Sigma_0$ is a fixed point of $f$ (at infinity).  The three  points of indeterminacy are obtained by finding the points of intersection of exceptional curves:
$$Ind(f) = \{[0:1:-1], e_1,e_2\}.$$
We see that $[0:1:-1]$ blows up to $\Sigma_0\ni[0:1:-1]$, so that this point will also be part of our intersection count when we are using the Lefschetz index.

Now let $\pi:Z\to {\bf P}^2$ denote the complex manifold obtained by blowing up $e_1$ and $e_2$.  If $H$ is the class of a generic hyperplane, and if $E_1$ and $E_2$ are the exceptional (blow-up) fibers, then $\{H,E_1,E_2\}$ is a basis for $H^2(Z)$.  Since $e_1$ and $e_2$ both blow up to lines, we find that $f_Z^*H=3H-E_1-E_2$.  Further, the preimages of $e_1$ and $e_2$ are lines which contain neither $e_1$ nor $e_2$.  Thus $f_Z^*E_1=f_Z^*E_2=H$.  

It follows that 
$$f^*_Z |H^*(Z)= \pmatrix{1 &&\cr
&\matrix{3&1&1\cr
-1&0&0\cr
-1&0&0\cr} & \cr
& & 2\cr}.$$

Now we consider the behavior of $f_Z:Z\to Z$.  The only exceptional hypersurfaces are $\{x_0=0\}$, which maps to a fixed point, and $\{Q(x)=0\}$, and there is only one point of indeterminacy: $ Ind(f_Z) = [0:1:-1]$.  For generic parameters, the orbit of $\{Q=0\}$ will be disjoint from $Ind(f_Z)$.  Thus, for generic parameters, we have $(f_Z^*)^n=(f_Z^n)^*$ for $n=1,2,3,\dots$  The number of  spurious periodic points for $f_Z$ is still 2, and so we subtract them from the index to obtain
$$\#Fix_n'=Tr(f_Z^{*n})  -2.$$

The characteristic polynomial of $f_Z^*|H^2(Z)$ is $x(x^2-3x+2)$, so $\tau_n=Tr(f_Z^{*n}|H^2)$ satisfies $\tau_{n+2}=3\tau_{n+1}-2\tau_n$, with initial condtions $\tau_1=3$ and $\tau_2=5$.  Thus $\tau_n=2^n+1$, and for generic parameters we have $\#Fix'_n = \tau_n+2^n-1=2^{n+1}$.
For arbitrary parameters, this gives an upper bound for $\#Fix'_n$.

\bigskip

\bigskip\centerline{\bf References}

\item{[BD]}  E. Bedford and J. Diller,  Real and complex dynamics of a family of birational maps of the plane: the golden mean subshift. Amer. J. Math. 127 (2005), no. 3, 595--646.

\item{[BK]}  E. Bedford and K-H Kim, Periodicities in linear fractional recurrences: degree growth of birational surface maps. Michigan Math. J. 54 (2006), no. 3, 647--670.

\item{[DF]} J. Diller and C. Favre, Dynamics of meromorphic maps of surfaces, Amer.\ J. of Math., 123 (2001), 1135--1169.

\item{[FS]} J.-E. Forn\ae ss and N. Sibony, Complex dynamics in higher dimension. II, pp.\  135--182, in {\sl Modern methods in complex analysis},  Ann.\ of Math.\ Stud.\  137, 1995.

\item{[Fu]}  W. Fulton, {\sl Intersection Theory}, Ergeb.\ Math.\ Grenzgeb., Springer-Verlag. 1984.
\item{[GL]}  E. Grove and G. Ladas, {\sl Periodicities in nonlinear difference equations.}  Advances in Discrete Mathematics and Applications, 4. Chapman \& Hall CRC, Boca Raton, FL, 2005.

\item{[KL]}  M. Kulenovi\'c and G. Ladas, {\sl Dynamics of second order rational difference equations. With open problems and conjectures.} Chapman \& Hall CRC, Boca Raton, FL, 2002.

\item{[KM]}   M. Kulenovi\'c and O. Merino, {\sl Discrete dynamical systems and difference equations with $Mathematica$.}  Chapman \& Hall CRC, Boca Raton, FL, 2002.

\bigskip
\rightline{Indiana University}

\rightline{Bloomington, IN 47405}

\rightline{\tt bedford@indiana.edu}

\medskip
\rightline{Florida State University}

\rightline{Tallahassee, FL 32306}

\rightline{\tt kim@math.fsu.edu}

\bye